\documentclass[graybox]{svmult}

\usepackage{type1cm}        
\usepackage{makeidx}         
\usepackage{graphicx}        
\usepackage{multicol}        
\usepackage[bottom]{footmisc}
\usepackage{newtxtext}       
\usepackage{newtxmath}       

\usepackage{enumerate} 
\usepackage[all]{xy}

\makeatletter
\newcommand\incircbin
{%
  \mathpalette\@incircbin
}
\newcommand\@incircbin[2]
{%
  \mathbin%
  {%
    \ooalign{\hidewidth$#1#2$\hidewidth\crcr$#1\bigcirc$}%
  }%
}

\makeatother

\makeindex

\begin{document}

\title*{The curvature of convex sums of metrics and applications}
\author{Leonardo F. Cavenaghi and Giovane Galindo and Llohann Speran\c{c}a}
\institute{Leonardo F. Cavenaghi \at Institute of Mathematics and Informatics, Bulgarian Academy of Sciences \email{leonardofcavenaghi@gmail.com}
\and Giovane Galindo \at Institute of Mathematics and Informatics, Bulgarian Academy of Sciences \email{giovanepgneto@gmail.com}
\and  Llohann Speran\c{c}a \at  \email{lsperanca@gmail.com}}

\maketitle

\abstract{In this note, we derive explicit formulae for the curvature of a convex sum of Riemannian metrics, \(g_t = (1-t)g_0 + t g_1\). We study whether such a deformation can increase the \emph{average} of the Riemann curvature component \(R_t(X,Y,Y,X)\) along an immersed, totally geodesic flat torus. Because a first-order increase is obstructed for a parallel tangent frame, we obtain necessary and sufficient conditions for \(g_t\) to have a positive average variation of order \(r \geq 2\). These conditions are applied to paths joining \(g_0\) to classical metric deformations, including conformal changes, vertical warpings, and Cheeger deformations.
}

\section{Introduction}\label{sec:1}

Let $(M,g_0)$ be a closed Riemannian manifold with non-negative sectional curvature containing a \emph{totally geodesic flat $2$-torus} $\mathbb{T}^2 \subset M$. In this note, we study the curvature behavior of metrics obtained from $g_0$ by a \emph{convex combination}
\[
g_t := (1-t)g_0 + t g_1, \qquad t \in [0,1],
\]
where $g_1$ is another metric on $M$.

Our work has two main goals. First, we derive a new explicit formula for the curvature tensor of $g_t$ by interpreting the convex sum as the metric induced on the graph of the identity in a product manifold (Section~\ref{sec:2})\footnote{Choi and Oh~\cite{choi2024injectivityradiuslowerbound} use the convex-sum curvature formula from an earlier version of this work in their study of curvature and injectivity-radius estimates for complete metrics.}. Second, we investigate whether such a deformation can increase, at some \emph{higher order}, the average of the Riemann curvature quantity $R_t(X,Y,Y,X)$ along the distinguished torus $\mathbb{T}^2$. More precisely, we seek conditions under which
\begin{equation}\label{eq:positiva}
\left.\frac{d^r}{dt^r}\right|_{t=0} \int_{\mathbb{T}^2} R_t(X,Y,Y,X)\, dA > 0, \qquad r \geq 2,
\end{equation}
where $X,Y$ form a parallel tangent frame on $\mathbb{T}^2$, and $dA$ denotes the area form induced by $g_0$. When \eqref{eq:positiva} holds, we say that $g_0$ admits an \emph{average positive $r$-order variation}.\\

This problem lies in the classical theory of curvature variation initiated by Bourguignon, Deschamps, and Sentenac \cite{Bourguignon1972, Bourguignon1973}. Their work reveals a fundamental tension: an analytic deformation starting at a compact Riemannian product metric without Killing vector fields cannot have positive sectional curvature for every sufficiently small $t>0$ \cite[Theorem 7.8]{Bourguignon1972}; on the other hand, under suitable technical assumptions, there exist polynomial deformations for which the first non-zero derivative of the curvature is positive on every mixed plane \cite[Theorem 4.1]{Bourguignon1973}. Our approach differs in two essential ways. We fix a specific deformation path, namely the convex combination $g_t$, and we ask for positivity only after averaging over a prescribed flat torus equipped with a parallel framing.\\

In section ~\ref{sec:3} we prove the following theorem that characterizes when such an average positive $r$-order variation occurs, giving a computable criterion in terms of the tensor $P$ defined by $g_0(P\cdot,\cdot)=g_1(\cdot,\cdot)$ and the difference of connections $\nabla^1-\nabla^0$.

\begin{theorem}\label{thm:rvariation-intro}
Assume that $X,Y$ span the tangent bundle of $\mathbb T^2$ and that $\nabla^0_X X = \nabla^0_X Y = \nabla^0_Y Y = 0$ on $\mathbb{T}^2$. Denote $R'_0(X,Y,Y,X) := \left.\frac{d}{dt}\right|_{t=0} R_t(X,Y,Y,X)$. Then
\begin{multline}\label{eq:first-variation-intro}
\int_{\mathbb{T}^2} \frac{1}{\|X\wedge Y\|_1^2} R'_0(X,Y,Y,X) \, dA_1 \\
= \int_{\mathbb{T}^2} \frac{1}{\|X\wedge Y\|_1^2} \Bigl[
g_1\bigl((\nabla^1_X X)^\top, (\nabla^1_Y Y)^\top\bigr) 
- g_1\bigl((\nabla^1_X Y)^\top, (\nabla^1_X Y)^\top\bigr) \Bigr] dA_1,
\end{multline}
where $dA_1$ is the area form induced on $\mathbb T^2$ by $g_1$ and $\|X\wedge Y\|_1^2 := g_1(X,X)g_1(Y,Y) - g_1(X,Y)^2$.

Moreover, if $R'_0(X,Y,Y,X) = 0$, then for any $r \geq 2$,
\begin{equation}\label{eq:higher-order-condition-intro}
\int_{\mathbb{T}^2} \left.\frac{d^r}{dt^r}\right|_{t=0} R_t(X,Y,Y,X) \, dA > 0 
\quad\Longleftrightarrow\quad 
\int_{\mathbb{T}^2} S_r^P(X,Y) \, dA < 0.
\end{equation}
Where

\begin{multline*}
S_r^P(X,Y) := g_1\!\bigl(P(I-P)^{r-2}(\nabla^1_X X - \nabla^0_X X), \nabla^1_Y Y - \nabla^0_Y Y\bigr) \\
- g_1\!\bigl(P(I-P)^{r-2}(\nabla^1_X Y - \nabla^0_X Y), \nabla^1_X Y - \nabla^0_X Y\bigr).
\end{multline*}
\end{theorem}

Finally, in Section~\ref{sec:4}, we apply it to three classical families of deformations:
\begin{itemize}
    \item \emph{Conformal changes} (Proposition~\ref{prop:conformal}). For $g_1=e^{2h}g_0$, suppose that $R'_0(X,Y,Y,X)=0$ pointwise on the torus. Then $h$ is constant there. The averaged variation of order $r\geq2$ is strictly positive exactly when this constant is positive, $r$ is odd, and the normal gradient of $h$ is nonzero somewhere on the torus. If the constant is zero, every averaged variation of order $r>2$ vanishes. Such changes cannot produce strictly positive curvature everywhere on a closed manifold.
    \item \emph{Vertical warpings} (Proposition~\ref{prop:gvw}). Under the same pointwise hypothesis, the basic warping function $f$ is constant on the torus. Every averaged variation of order $r\geq2$ therefore vanishes there, although vertical warping may introduce negative curvature elsewhere \cite{speranca_oddbundles}.
    \item \emph{Cheeger deformations} (Proposition~\ref{prop:CD}), for which we analyze the limiting behavior for actions locally free along the torus as the deformation parameter tends to infinity.
\end{itemize}

Although condition \eqref{eq:positiva} provides a concrete mechanism for increasing a key curvature quantity on a flat torus, the construction of globally positively curved metrics remains a subtle problem with well-known obstructions \cite{wilking2002manifolds}. The results presented here provide precise tools for controlling curvature under deformation, even though additional global input is still needed to produce new examples of positively curved manifolds.

\section{The shape operator of a convex sum of Riemannian metrics}\label{sec:2}

In this section we will derive a new formula for the curvature of a convex sum of metrics.\\

Let $(M,g_0)$ and $(N,g_N)$ be Riemannian manifolds of dimension $m$ and $n$ respectively, and let $r : M \to N$ be a smooth map. The graph of $r$,
\[
\Gamma_r := \{(p,r(p)) \in M\times N : p \in M\},
\]
is a submanifold of $M\times N$. The map $\Xi : M \to \Gamma_r$ defined by $p \mapsto (p,r(p))$ is an embedding identifying $M$ with $\Gamma_r$.

Now specialize to the case $M = N$ (so $m = n$), $g_N = g_1$, and $r = \operatorname{Id}_M$. Let $f : M\to (0,1)$ be a smooth function and consider the two metrics $(1-f)g_0$ and $fg_1$ on $M$. Equip the product $M\times M$ with the Riemannian metric
\[
\overline{g} := \pi_1^*\big((1-f)g_0\big) + \pi_2^*\big(fg_1\big),
\]
where $\pi_1,\pi_2 : M\times M \to M$ are the projections onto the first and second factors respectively. Then $\Gamma_{\operatorname{Id}_M}$ is a submanifold of $(M\times M,\overline{g})$, and the induced metric on $\Gamma_{\operatorname{Id}_M}$ via $\Xi$ is
\begin{align*}
\Xi^*\overline{g} 
&= (\pi_1\circ\Xi)^*\big((1-f)g_0\big) + (\pi_2\circ\Xi)^*\big(fg_1\big) \\
&= (1-f)g_0 + fg_1,
\end{align*}
because $\pi_1\circ\Xi = \pi_2\circ\Xi = \operatorname{Id}_M$. Thus the convex sum $(1-f)g_0 + fg_1$ coincides with the metric induced on $M$ by the embedding $\Xi$. Consequently, understanding the curvature of the convex sum reduces to computing the second fundamental form of $\Xi(M) \subset (M\times M,\overline{g})$.

Since $g_0$ and $g_1$ are Riemannian metrics on $M$, at each point $p\in M$ there exists a unique symmetric, positive-definite endomorphism $P : T_pM \to T_pM$ satisfying
\begin{equation}\label{eq:P-definition}
(1-f)g_0(PX,Y) = fg_1(X,Y) \qquad \forall\, X,Y\in T_pM.
\end{equation}
We use this tensor $P$ to encode the difference of the two metrics and we will use it to describe the curvature of $g_t$.
\begin{lemma}\label{lem:tangent-normal}
The tangent and normal spaces of $\Gamma_{\operatorname{Id}_M}$ at $(p,p)$ are given by
\begin{align*}
T_{(p,p)}\Gamma_{\operatorname{Id}_M} &= \{(X,X) : X\in T_pM\}, \\
\nu_{(p,p)}\Gamma_{\operatorname{Id}_M} &= \{(-PX,X) : X\in T_pM\},
\end{align*}
with respect to the metric $\overline{g}$. Moreover, the inverse of the isomorphism
\[
\Phi : T_pM\times T_pM \to T_{(p,p)}(M\times M), \quad (X,Y) \mapsto (X-PY,\, X+Y)
\]
is given by
\begin{equation}\label{eq:inverse-iso}
\Phi^{-1}\begin{pmatrix}U\\V\end{pmatrix} = 
\begin{pmatrix}
(1+P)^{-1} & P(1+P)^{-1} \\[2pt]
-(1+P)^{-1} & (1+P)^{-1}
\end{pmatrix}
\begin{pmatrix}U\\V\end{pmatrix}.
\end{equation}
\end{lemma}

\begin{proof}
That $T_{(p,p)}\Gamma_{\operatorname{Id}_M} = \{(X,X) : X\in T_pM\}$ follows because $\Xi(p)=(p,p)$ and $d\Xi(X)=(X,X)$. For orthogonality, take $(X,X)\in T_{(p,p)}\Gamma_{\operatorname{Id}_M}$ and $(-PY,Y)\in \nu_{(p,p)}\Gamma_{\operatorname{Id}_M}$:
\begin{align*}
\overline{g}\big((X,X),(-PY,Y)\big) 
&= (1-f)g_0(X,-PY) + fg_1(X,Y) \\
&= -fg_1(X,Y) + fg_1(X,Y) = 0,
\end{align*}
using \eqref{eq:P-definition}. Since $\dim\nu_{(p,p)}\Gamma_{\operatorname{Id}_M} = m$, the claim follows.

To obtain \eqref{eq:inverse-iso}, solve $(U,V) = (X-PY,\, X+Y)$ for $X$ and $Y$:
\[
\begin{cases}
U = X - PY, \\
V = X + Y.
\end{cases}
\]
Subtracting gives $V-U = (1+P)Y$, hence $Y = (1+P)^{-1}(V-U)$. Substituting back yields $X = V - Y = (1+P)^{-1}(U+PV)$.
\end{proof}

Set $O := (1+P)^{-1}$.

\begin{lemma}\label{lem:normal-projection}
The orthogonal projection onto $\nu\Gamma_{\operatorname{Id}_M}$ is given by
\begin{equation}\label{eq:normal-projection}
\operatorname{proj}_{\nu\Gamma_{\operatorname{Id}_M}}\begin{pmatrix}X\\Y\end{pmatrix} = 
\begin{pmatrix}
P O X - P O Y \\[2pt]
-O X + O Y
\end{pmatrix}.
\end{equation}
\end{lemma}

\begin{proof}
Using Lemma~\ref{lem:tangent-normal}:
\[
\operatorname{proj}_{\nu\Gamma_{\operatorname{Id}_M}} = \chi'\circ\operatorname{proj}_{(T_pM,fg_1)}\circ\Phi^{-1},
\]
where $\chi' : (T_pM,fg_1) \to \nu_{(p,p)}\Gamma_{\operatorname{Id}_M}$ is given by $Y \mapsto (-PY, Y)$, and $\operatorname{proj}_{(T_pM,fg_1)}$ denotes projection onto the second factor. A matrix multiplication yields \eqref{eq:normal-projection}.
\end{proof}

\begin{lemma}[Shape operator of the convex sum]\label{lem:shape-operator}
Let $\sigma$ denote the second fundamental form of $\Gamma_{\operatorname{Id}_M}\subset(M\times M,\overline{g})$. Define the difference tensor
\[
D(X,Y) := \nabla'_X Y - \nabla_X Y,
\]
where $\nabla'$ and $\nabla$ are the Levi-Civita connections of $fg_1$ and $(1-f)g_0$, respectively. Then for any $X,Y,X',Y'\in T_pM$,
\begin{multline}\label{eq:shape-formula}
\overline{g}\Bigl(\sigma\bigl(d\Xi(X),d\Xi(Y)\bigr),\,
\sigma\bigl(d\Xi(X'),d\Xi(Y')\bigr)\Bigr) \\
= fg_1\Bigl((1+P)^{-1}D(X,Y),\, D(X',Y')\Bigr).
\end{multline}
\end{lemma}

\begin{proof}
By definition,
\[
\sigma(d\Xi(X),d\Xi(Y)) = \operatorname{proj}_{\nu\Gamma_{\operatorname{Id}_M}}\bigl(\nabla_X Y,\; \nabla'_X Y\bigr).
\]
Using Lemma~\ref{lem:normal-projection},
\begin{align*}
\sigma(d\Xi(X),d\Xi(Y)) 
&= \begin{pmatrix}
P O \nabla_X Y - P O \nabla'_X Y \\[2pt]
-O \nabla_X Y + O \nabla'_X Y
\end{pmatrix} \\
&= \begin{pmatrix}
-P O D(X,Y) \\[2pt]
O D(X,Y)
\end{pmatrix}.
\end{align*}
Similarly, $\sigma(d\Xi(X'),d\Xi(Y')) = \begin{pmatrix} -P O D(X',Y') \\ O D(X',Y') \end{pmatrix}$. Their inner product is
\begin{align*}
&\overline{g}\Bigl(\sigma(d\Xi(X),d\Xi(Y)),\sigma(d\Xi(X'),d\Xi(Y'))\Bigr) \\
&= (1-f)g_0\bigl(P O D(X,Y), P O D(X',Y')\bigr) 
   + fg_1\bigl(O D(X,Y), O D(X',Y')\bigr) \\
&= fg_1\bigl(O D(X,Y), P O D(X',Y')\bigr) 
   + fg_1\bigl(O D(X,Y), O D(X',Y')\bigr) \quad\text{(by \eqref{eq:P-definition})} \\
&= fg_1\bigl(O D(X,Y), (1+P)O D(X',Y')\bigr) \\
&= fg_1\bigl(D(X,Y), O D(X',Y')\bigr),
\end{align*}
since $(1+P)O = I$ and $O$ is symmetric with respect to $fg_1$.
\end{proof}

\section{Any order variations of convex sums of metrics} \label{sec:3}

Now we apply the results of last section to obtain necessary and sufficient conditions for an average positive $r$-order variation to occur in convex sums.\\

Let $(M,g_0)$ be a closed Riemannian manifold containing a totally geodesic flat $2$-torus $\mathbb{T}^2 \subset M$. We denote by $X,Y$ tangent vector fields to $\mathbb{T}^2$. Let $g_1$ be another metric on $M$ and consider the convex combination $g_t := (1-t)g_0 + t g_1$ for $t \in [0,1]$.

\begin{lemma}\label{lem:resolvetudo}
Assume that $\nabla^0_X X = \nabla^0_X Y = \nabla^0_Y Y = 0$ on $\mathbb{T}^2$. Denote by $\nabla^t$ the Levi-Civita connection of $g_t$. Then for sufficiently small $t$ (such that $\|t(I-P)\| < 1$), the Riemann curvature component admits the expansion
\begin{align}
R_t(X,Y,Y,X) =& (1-t)R_0(X,Y,Y,X) + tR_1(X,Y,Y,X) \nonumber \\
&+ t\Bigl[g_1\bigl((\nabla^1_X X)^\perp, (\nabla^1_Y Y)^\perp\bigr) 
     - g_1\bigl((\nabla^1_X Y)^\perp, (\nabla^1_X Y)^\perp\bigr)\Bigr] \nonumber \\
&+ t\Bigl[g_1\bigl((\nabla^1_X X)^\top, (\nabla^1_Y Y)^\top\bigr) 
     - g_1\bigl((\nabla^1_X Y)^\top, (\nabla^1_X Y)^\top\bigr)\Bigr] \nonumber \\
&- t^2\sum_{n=1}^\infty \Bigl[g_1\bigl(P(t(I-P))^{n-1}\nabla^1_X X, \nabla^1_Y Y\bigr) \nonumber \\
&\qquad\qquad - g_1\bigl(P(t(I-P))^{n-1}\nabla^1_X Y, \nabla^1_X Y\bigr)\Bigr],
\end{align}
where $^\top$ and $^\perp$ denote the $g_1$-orthogonal projections onto $T\mathbb{T}^2$ and its orthogonal complement, respectively, and $P$ is defined by $g_0(P\cdot,\cdot)=g_1(\cdot,\cdot)$.
\end{lemma}

\begin{proof}
The starting point is the expression for $R_t(X,Y,Y,X)$ from the Gauss equation for the graph embedding $\Xi: (M,g_t) \hookrightarrow (M\times M, (1-t)g_0 \times t g_1)$. For $0<t<1$, from Lemma~\ref{lem:shape-operator} we have:
\begin{align*}
R_t(X,Y,Y,X) &= (1-t)R_0(X,Y,Y,X) + tR_1(X,Y,Y,X) \\
 &\quad - g_1(tO_t D(X,Y),D(Y,X)) + g_1(tO_t D(X,X),D(Y,Y)),
\end{align*}
where $D(U,V) = \nabla^1_U V - \nabla^0_U V$ and $O_t = (I + \frac{t}{1-t}P)^{-1}$.

For $(U,V)=(X,X),(X,Y),(Y,Y)$, the hypotheses give $D(U,V)=\nabla^1_U V$. The operator $O_t$ can be expanded as a power series. First note that:
\[
O_t = \left(I + \frac{t}{1-t}P\right)^{-1} = (1-t)\left(I - t(I-P)\right)^{-1} = (1-t)\sum_{n=0}^\infty t^n (I-P)^n,
\]
provided $\|t(I-P)\| < 1$. Then:
\begin{align*}
t O_t &= t(1-t)\sum_{n=0}^\infty t^n (I-P)^n \\
&= t \sum_{n=0}^\infty t^n(I-P)^n - t^2 \sum_{n=0}^\infty t^n(I-P)^n \\
&= t + \sum_{n=1}^\infty t^{n+1}(I-P)^n - \sum_{n=1}^\infty t^{n+1}(I-P)^{n-1} \\
&= t + \sum_{n=1}^\infty t^{n+1} \left( (I-P)^n - (I-P)^{n-1} \right) \\
&= t - P\sum_{n=1}^\infty t^{n+1} (I-P)^{n-1}.
\end{align*}

Writing $A_t=tO_t$, substituting this into the curvature expression yields:
\begin{align*}
R_t(X,Y,Y,X) &= (1-t)R_0(X,Y,Y,X) + tR_1(X,Y,Y,X) \\
&\quad +g_1(A_t\nabla^1_X X,\nabla^1_Y Y)-g_1(A_t\nabla^1_X Y,\nabla^1_X Y).
\end{align*}
Decomposing $\nabla^1_U V$ into tangential and normal components and using orthogonality completes the proof.
\end{proof}

\begin{remark}
    Note that, because $\mathbb{T}^2$ is a totally geodesic flat torus, a tangent frame $X,Y$ satisfying $\nabla^0_X X = \nabla^0_X Y = \nabla^0_Y Y = 0$ always exists.
\end{remark}

Define for $r \geq 2$ the quantity
\begin{multline*}
S_r^P(X,Y) := g_1\!\bigl(P(I-P)^{r-2}(\nabla^1_X X - \nabla^0_X X), \nabla^1_Y Y - \nabla^0_Y Y\bigr) \\
- g_1\!\bigl(P(I-P)^{r-2}(\nabla^1_X Y - \nabla^0_X Y), \nabla^1_X Y - \nabla^0_X Y\bigr).
\end{multline*}

\begin{theorem}\label{thm:rvariation}
Assume that $X,Y$ span the tangent bundle of $\mathbb T^2$ and that $\nabla^0_X X = \nabla^0_X Y = \nabla^0_Y Y = 0$ on $\mathbb{T}^2$. Denote $R'_0(X,Y,Y,X) := \left.\frac{d}{dt}\right|_{t=0} R_t(X,Y,Y,X)$. Then
\begin{multline}\label{eq:first-variation}
\int_{\mathbb{T}^2} \frac{1}{\|X\wedge Y\|_1^2} R'_0(X,Y,Y,X) \, dA_1 \\
= \int_{\mathbb{T}^2} \frac{1}{\|X\wedge Y\|_1^2} \Bigl[
g_1\bigl((\nabla^1_X X)^\top, (\nabla^1_Y Y)^\top\bigr) 
- g_1\bigl((\nabla^1_X Y)^\top, (\nabla^1_X Y)^\top\bigr) \Bigr] dA_1,
\end{multline}
where $dA_1$ is the area form induced on $\mathbb T^2$ by $g_1$ and $\|X\wedge Y\|_1^2 := g_1(X,X)g_1(Y,Y) - g_1(X,Y)^2$.

Moreover, if $R'_0(X,Y,Y,X) = 0$, then for any $r \geq 2$,
\begin{equation}\label{eq:higher-order-condition}
\int_{\mathbb{T}^2} \left.\frac{d^r}{dt^r}\right|_{t=0} R_t(X,Y,Y,X) \, dA > 0 
\quad\Longleftrightarrow\quad 
\int_{\mathbb{T}^2} S_r^P(X,Y) \, dA < 0.
\end{equation}
\end{theorem}

\begin{proof}
Differentiating the expansion in Lemma~\ref{lem:resolvetudo} at $t=0$ yields:
\begin{align*}
R'_0(X,Y,Y,X) &= R_1(X,Y,Y,X) - R_0(X,Y,Y,X) \\
&\quad + g_1\bigl((\nabla^1_X X)^\perp, (\nabla^1_Y Y)^\perp\bigr) 
   - g_1\bigl((\nabla^1_X Y)^\perp, (\nabla^1_X Y)^\perp\bigr) \\
&\quad + g_1\bigl((\nabla^1_X X)^\top, (\nabla^1_Y Y)^\top\bigr) 
   - g_1\bigl((\nabla^1_X Y)^\top, (\nabla^1_X Y)^\top\bigr).
\end{align*}
Since $\mathbb{T}^2$ is flat with respect to $g_0$, we have $R_0(X,Y,Y,X)=0$. Moreover, by the Gauss equation for the submanifold $\mathbb{T}^2 \subset M$ with the metric $g_1$, we have:
\[
R_1(X,Y,Y,X) = K^1_{\mathbb{T}^2} \|X\wedge Y\|_1^2 - \left[g_1((\nabla^1_X X)^\perp, (\nabla^1_Y Y)^\perp) - g_1((\nabla^1_X Y)^\perp, (\nabla^1_X Y)^\perp)\right],
\]
where $K^1_{\mathbb{T}^2}$ is the intrinsic Gaussian curvature of $\mathbb{T}^2$ with respect to $g_1$. Substituting this gives:
\[
R'_0(X,Y,Y,X) = K^1_{\mathbb{T}^2} \|X\wedge Y\|_1^2 + g_1((\nabla^1_X X)^\top, (\nabla^1_Y Y)^\top) - g_1((\nabla^1_X Y)^\top, (\nabla^1_X Y)^\top).
\]
Dividing by $\|X\wedge Y\|_1^2$, integrating over $\mathbb{T}^2$ and applying the Gauss-Bonnet theorem ($\int_{\mathbb{T}^2} K^1_{\mathbb{T}^2} \, dA_1 = 0$) yields \eqref{eq:first-variation}.

Now assume $R'_0(X,Y,Y,X)=0$. From Lemma~\ref{lem:resolvetudo}, for $r \geq 2$, the $r$th derivative at $t=0$ of the sum $ -t^2\sum_{n=1}^\infty \left[g_1(P(t(I-P))^{n-1}\nabla^1_X X, \nabla^1_Y Y) - g_1(P(t(I-P))^{n-1}\nabla^1_X Y, \nabla^1_X Y)\right] $ receives contributions only from the term where the power of $t$ is exactly $r$. This occurs when $n+1 = r$, i.e., $n = r-1$. For this term,
\begin{align*}
&\left.\frac{d^r}{dt^r}\right|_{t=0} \left( -t^2 \cdot \left[g_1(P(t(I-P))^{r-2}\nabla^1_X X, \nabla^1_Y Y) - \cdots \right] \right) \\
&= \left.\frac{d^r}{dt^r}\right|_{t=0} \left( -t^{r} \left[g_1(P(I-P)^{r-2}\nabla^1_X X, \nabla^1_Y Y) - \cdots \right] \right) \\
&= -r! \left[g_1(P(I-P)^{r-2}\nabla^1_X X, \nabla^1_Y Y) - g_1(P(I-P)^{r-2}\nabla^1_X Y, \nabla^1_X Y)\right].
\end{align*}
This simplifies to:
\[
\left.\frac{d^r}{dt^r}\right|_{t=0} R_t(X,Y,Y,X) = -r! \cdot S_r^P(X,Y).
\]
Integrating over $\mathbb{T}^2$ yields the equivalence \eqref{eq:higher-order-condition}.
\end{proof}

\section{Curvature variations along paths joining a metric to classical deformations}\label{sec:4}

Throughout this section, we maintain the assumptions of previous sections. We consider a closed Riemannian manifold $(M,g_0)$ with a totally geodesic flat torus $\mathbb{T}^2 \subset M$. We always denote by $\{X,Y\}$ a $g_0$-parallel orthonormal frame tangent to $\mathbb{T}^2$. In addition, we shall assume that $\nabla^0_XY = \nabla^0_XX = \nabla^0_YY = 0$ and that $R'_0(X,Y,Y,X)=0$.

\subsection{Conformal changes}

\begin{proposition}[Conformal class]\label{prop:conformal}
Let $g_1 = f g_0 = e^{2h} g_0$, where $h : M \to \mathbb{R}$ is a smooth function and $f = e^{2h}$. Then $f|_{\mathbb{T}^2}=c$ is constant and for all $r \geq 2$,
\begin{equation*}
\int_{\mathbb{T}^2}\left.\frac{d^r}{dt^r}\right|_{t=0} R_t(X,Y,Y,X)\,dA
= -r!c^2(1-c)^{r-2}\int_{\mathbb{T}^2}\|(\nabla^0h)^\perp\|_{g_0}^2\,dA .
\end{equation*}
The averaged $r$-th variation is strictly positive if and only if $h|_{\mathbb T^2}>0$, $r$ is odd, and $(\nabla^0h)^\perp$ is nonzero somewhere on $\mathbb T^2$. If $h|_{\mathbb T^2}=0$, every variation of order $r>2$ vanishes.
\end{proposition}

\begin{proof}
For a conformal change $g_1 = f g_0$, the tensor $P$ defined by $g_0(P\cdot,\cdot)=g_1(\cdot,\cdot)$ is simply multiplication by $f$, i.e., $P = f \operatorname{Id}$. The Levi-Civita connection $\nabla^f$ of $g_1 = f g_0$ is related to $\nabla^0$ by (see \cite[Theorem 1.159]{besse})
\[
\nabla^f_U V = \nabla^0_U V + dh(U)V + dh(V)U - g_0(U,V)\nabla^0 h.
\]

Using that $\nabla^0_X X = \nabla^0_X Y = \nabla^0_Y Y = 0$, a direct computation yields
\begin{align}
g_0(\nabla^f_X X, \nabla^f_Y Y) - \|\nabla^f_X Y\|_{g_0}^2 
&= -2\|(\nabla^0 h)^\top\|_{g_0}^2 + \|(\nabla^0 h)^\perp\|_{g_0}^2, \label{eq:conformal-computation}
\end{align}
where $(\nabla^0 h)^\top$ and $(\nabla^0 h)^\perp$ denote the components of $\nabla^0 h$ tangent and orthogonal to $\mathbb{T}^2$, respectively.

Since $R'_0(X,Y,Y,X) = 0$, Theorem~\ref{thm:rvariation} implies via equation \eqref{eq:first-variation} that
\[
\int_{\mathbb{T}^2} \frac{1}{\|X\wedge Y\|_1^2} \Bigl[g_1\bigl((\nabla^f_X X)^\top, (\nabla^f_Y Y)^\top\bigr) - g_1\bigl((\nabla^f_X Y)^\top, (\nabla^f_X Y)^\top\bigr) \Bigr] dA_1 = 0.
\]
Given the conformal connection formula, this forces $\int_{\mathbb{T}^2} \|(\nabla^0 h)^\top\|_{g_0}^2 \, dA = 0$, and hence $(\nabla^0 h)^\top = 0$ on $\mathbb{T}^2$, meaning $h$ (and therefore $f$) is constant on $\mathbb{T}^2$. Let $c := f|_{\mathbb{T}^2} > 0$.

Now, by Theorem~\ref{thm:rvariation}, for $r \geq 2$,
\begin{align*}
\int_{\mathbb{T}^2} \left.\frac{d^r}{dt^r}\right|_{t=0} R_t(X,Y,Y,X) \, dA
&= -r!\int_{\mathbb{T}^2} S_r^P(X,Y) \, dA \\
&= -r!\int_{\mathbb{T}^2} f^2(1-f)^{r-2} \bigl[g_0(\nabla^f_X X, \nabla^f_Y Y) - \|\nabla^f_X Y\|_{g_0}^2\bigr] dA \\
&= -r!c^2(1-c)^{r-2} \int_{\mathbb{T}^2} \left(-2\|(\nabla^0 h)^\top\|_{g_0}^2 + \|(\nabla^0 h)^\perp\|_{g_0}^2\right) dA \\
&= -r! c^2(1-c)^{r-2} \int_{\mathbb{T}^2} \|(\nabla^0 h)^\perp\|_{g_0}^2 \, dA.
\end{align*}
This concludes the proof.
\end{proof}

\subsection{Vertical warping}

Now assume $(M,g_0)$ is equipped with a Riemannian foliation $\mathcal{F}$. Denote by $\mathcal{V}$ the tangent bundle to the leaves and by $\mathcal{H}$ its $g_0$-orthogonal complement. For a vector $U$, write $U = U^{\mathcal{V}} + U^{\mathcal{H}}$ for its vertical and horizontal components. Let $\sigma : \mathcal{V} \times \mathcal{V} \to \mathcal{H}$ be the shape operator of the leaves, and $A : \mathcal{H} \times \mathcal{H} \to \mathcal{V}$ the O'Neill tensor, defined by $A_{X^{\mathcal{H}}} Y^{\mathcal{H}} = \frac{1}{2}[X^{\mathcal{H}}, Y^{\mathcal{H}}]^{\mathcal{V}}$.

\begin{proposition}[General vertical warping]\label{prop:gvw}
Assume the flat torus $\mathbb{T}^2$ is spanned at each point by a horizontal vector $X$ and a vertical vector $Y$. Let $f : M \to \mathbb{R}$ be a basic function (constant on leaves) and consider the vertically warped metric
\[
g_f := g_0|_{\mathcal{H}} + e^{2f} g_0|_{\mathcal{V}}.
\]
Under the standing hypothesis $R'_0(X,Y,Y,X)=0$ on $\mathbb T^2$, the restriction of $f$ to the torus is constant. In particular, there is no positive averaged variation of order $r\geq2$, since
\[
\int_{\mathbb{T}^2} \left.\frac{d^r}{dt^r}\right|_{t=0} R_t(X,Y,Y,X) \, dA=0.
\]
\end{proposition}

\begin{proof}
For the warped metric $g_f$, the tensor $P$ satisfies $P|_{\mathcal{H}} = \operatorname{Id}$ and $P|_{\mathcal{V}} = e^{2f} \operatorname{Id}$. The Levi-Civita connection $\nabla^f$ of $g_f$ is given by (see \cite[Chapter 2]{gw})
\begin{align}
\nabla^f_{X^{\mathcal{H}}} Y^{\mathcal{H}} &= \nabla^0_{X^{\mathcal{H}}} Y^{\mathcal{H}}, \label{eq:gw1} \\
\nabla^f_{X^{\mathcal{V}}} Y^{\mathcal{H}} &= \nabla^0_{X^{\mathcal{V}}} Y^{\mathcal{H}} + (1 - e^{2f}) A^*_{Y^{\mathcal{H}}} X^{\mathcal{V}} + df(Y^{\mathcal{H}}) X^{\mathcal{V}}, \label{eq:gw2} \\
\nabla^f_{X^{\mathcal{H}}} Y^{\mathcal{V}} &= \nabla^0_{X^{\mathcal{H}}} Y^{\mathcal{V}} + (1 - e^{2f}) A^*_{X^{\mathcal{H}}} Y^{\mathcal{V}} + df(X^{\mathcal{H}}) Y^{\mathcal{V}}, \label{eq:gw3} \\
\nabla^f_{X^{\mathcal{V}}} Y^{\mathcal{V}} &= \nabla^0_{X^{\mathcal{V}}} Y^{\mathcal{V}} + (e^{2f} - 1) \sigma(X^{\mathcal{V}}, Y^{\mathcal{V}}) - e^{2f} g_0(X^{\mathcal{V}}, Y^{\mathcal{V}}) \nabla^0 f. \label{eq:gw4}
\end{align}

Since $X = X^{\mathcal{H}}$ and $Y = Y^{\mathcal{V}}$ are orthonormal with respect to $g_0$, and $\nabla^0_X X = \nabla^0_X Y = \nabla^0_Y Y = 0$, the horizontal part of $\nabla^0_XY=0$ gives $A_X^*Y=0$. Hence formulas \eqref{eq:gw1}--\eqref{eq:gw4} yield
\begin{align*}
\nabla^f_X X &= 0, \\
\nabla^f_Y Y &= (e^{2f} - 1) \sigma(Y,Y) - e^{2f} \nabla^0 f, \\
\nabla^f_X Y &= df(X) Y.
\end{align*}

Now compute $S_r^P(X,Y)$ as defined before:
\begin{align*}
S_r^P(X,Y) &= g_1\!\left(P(I-P)^{r-2} \nabla^f_X X, \nabla^f_Y Y\right) 
          - g_1\!\left(P(I-P)^{r-2} \nabla^f_X Y, \nabla^f_X Y\right) \\
&= e^{2f} g_0\!\left(e^{2f}(1 - e^{2f})^{r-2} \nabla^f_X X, \nabla^f_Y Y\right) \\
&\quad - e^{2f} g_0\!\left(e^{2f}(1 - e^{2f})^{r-2} \nabla^f_X Y, \nabla^f_X Y\right) \\
&= -e^{4f}(1 - e^{2f})^{r-2} \|df(X) Y\|_{g_0}^2.
\end{align*}

Write $F=e^{2f}$. In local parallel coordinates with $X=\partial_x$ and $Y=\partial_y$, the induced metric on the torus is $dx^2+(1+t(F-1))\,dy^2$. Its first curvature variation is $-X^2F/2$. By the Gauss equation and total geodesy at $t=0$, this equals $R'_0(X,Y,Y,X)$, which vanishes by hypothesis. Thus $F$ is affine along each complete $X$-trajectory. Since $F$ is bounded on the torus, $XF=0$, hence $Xf=0$. Basicness gives $Yf=0$, so $f$ is constant on the torus. The expression above for $S_r^P(X,Y)$ now vanishes pointwise, and the curvature expansion in Lemma~\ref{lem:resolvetudo} gives the asserted vanishing for every $r\geq2$.
\end{proof}

\subsection{Cheeger deformations}

Let $(M,g_0)$ be a Riemannian manifold on which a compact Lie group $G$ acts by isometries. Fix a bi-invariant metric $Q$ on $\mathfrak{g}$, the Lie algebra of $G$. For each $p \in M$, let $\mathfrak{g}_p$ be the kernel of the differential of the orbit map; its $Q$-orthogonal complement $\mathfrak{m}_p$ is identified with the tangent space to the orbit $G\cdot p$. The \emph{orbit tensor} $\mathcal{O}_p : \mathfrak{m}_p \to \mathfrak{m}_p$ is defined by
\[
g_0(V^*, W^*) = Q(\mathcal{O}_p V, W), \quad V,W \in \mathfrak{m}_p,
\]
where $V^*$ denotes the action field of $V$.

For $s \geq 0$, the \emph{Cheeger deformation} of $g_0$ is the metric $g_s$ defined by
\[
g_s(X,Y) = g_0(P_s X, Y), \quad \text{where } P_s|_{\mathcal{H}} = \operatorname{Id}, \quad P_s|_{\mathcal{V}} = (1 + s \mathcal{O})^{-1}.
\]
Here $\mathcal{V}$ is the vertical distribution (tangent to orbits) and $\mathcal{H}$ is the horizontal distribution ($g_0$-orthogonal to orbits).

Consider the two-parameter family $g_t^s := (1-t)g_0 + t g_s$ and denote by $R_t^s$ the curvature tensor of $g_t^s$.

\begin{proposition}[Cheeger deformations]\label{prop:CD}
Let $X,Y$ be tangent to $\mathbb{T}^2$ as before, with $X = X^{\mathcal{H}}$, $Y = Y^{\mathcal{V}}$. Assume that the action is locally free along $\mathbb T^2$. Write $Y=y^*$ pointwise and set $b=(d\mathcal O)(X)\mathcal O y$, differentiating the orbit tensor in a fixed $Q$-orthonormal basis of $\mathfrak g$. For $r \geq 3$, define the rescaled vectors $\tilde{X}_s := P_s^{-1}X$, $\tilde{Y}_s := P_s^{-1}Y$. Then
\begin{equation*}
\lim_{s\to\infty} \int_{\mathbb{T}^2} \left.\frac{d^r}{dt^r}\right|_{t=0} R_t^s(\tilde{X}_s, \tilde{Y}_s, \tilde{Y}_s, \tilde{X}_s) \, dA > 0 
\Longleftrightarrow \quad b\ne0 \text{ somewhere on }\mathbb{T}^2 .
\end{equation*}
\end{proposition}

\begin{proof}
In the chosen basis put $H=\mathcal O$, $H'=(dH)(X)$ and $L_s=(I+sH)^{-1}$. The vertical metric of $g_s$ has matrix $H_s=HL_s=(H^{-1}+sI)^{-1}$. Write $D^s=\nabla^{g_s}-\nabla^0$. A horizontal $g_0$-geodesic remains a $g_s$-geodesic by the Cheeger submersion, so $D^s(X,X)=0$.

For fixed action fields $u^*,z^*$, the Killing equations give
$2g_0(\nabla^0_Xu^*,z^*)=Xg_0(u^*,z^*)$, since $[u^*,z^*]$ is vertical. The same identity holds for $g_s$. Differentiating $H_s^{-1}=H^{-1}+sI$ therefore gives, in action coordinates,
\[
(D^s(X,u^*))^{\mathcal V}
 =\tfrac12(H_s^{-1}H_s'-H^{-1}H')u
 =\tfrac12 H^{-1}H'(L_s-I)u .
\]
For $\tilde Y_s$ take $u=L_s^{-1}y$ at the point of evaluation. The last expression becomes $-\tfrac{s}{2}H^{-1}b$. Also $\tilde X_s=X$. For $r\ge3$, $P_s(I-P_s)^{r-2}$ kills the horizontal component. Applying the diagonal Gauss formula in the proof of Lemma~\ref{lem:resolvetudo} to these fixed vectors and using $sHL_s=I-L_s$ gives
\[
\left.\frac{d^r}{dt^r}\right|_{t=0} R_t^s(\tilde X_s,\tilde Y_s,\tilde Y_s,\tilde X_s)
=\frac{r!}{4}\,Q\big(H^{-3}(I-L_s)^r b,b\big).
\]
Local freeness and compactness bound the least eigenvalue of $H$ away from zero on the torus. Thus $L_s\to0$ uniformly and the integral converges to $\tfrac{r!}{4}\int_{\mathbb T^2}Q(H^{-3}b,b)\,dA$. Positive definiteness of $H^{-3}$ gives the equivalence.
\end{proof}

\section*{Acknowledgments}
L.F.Cavenaghi and G.Galindo are supported by the Simons Foundation, grant SFI-MPS-T-Institutes-00007697, and the Ministry of Education and Science of the Republic of Bulgaria, grant DO1-239/10.12.2024 .

L. D. Sperança has no current academic affiliation.

\bibliographystyle{abbrv}
 \bibliography{main}

@article{choi2024injectivityradiuslowerbound,
title = {Injectivity radius lower bound of convex sum of tame Riemannian metrics and applications to symplectic topology},
journal = {Advances in Mathematics},
volume = {479},
pages = {110443},
year = {2025},
issn = {0001-8708},
doi = {https://doi.org/10.1016/j.aim.2025.110443},
url = {https://www.sciencedirect.com/science/article/pii/S000187082500341X},
author = {Jaeyoung Choi and Yong-Geun Oh}
}

@article{Bourguignon1972,
  author    = {Bourguignon, J.-P. and Deschamps, A. and Sentenac, P.},
  title     = {{Conjecture de H. Hopf sur les produits de variétés}},
  journal   = {Annales scientifiques de l'École Normale Supérieure},
  series    = {Série 4},
  volume    = {5},
  number    = {2},
  pages     = {277--302},
  year      = {1972},
  doi       = {10.24033/asens.1229},
  url       = {http://www.numdam.org/articles/10.24033/asens.1229/}
}

@article{Bourguignon1973,
author = {Bourguignon, J.-P. and  Deschamps, A. and Sentenac, P.},
journal = {Annales scientifiques de l'École Normale Supérieure},
language = {fre},
number = {1},
pages = {1-16},
publisher = {Elsevier},
title = {Quelques variations particulières d'un produit de métriques},
url = {http://eudml.org/doc/81911},
volume = {6},
year = {1973},
}

@book{besse,
	author = {Besse, A.},
	title = {Einstein Manifolds},
	year = {2008},
	publisher ={Springer}
}

@Book{gw,
	author = "Gromoll, D. and Walschap, G.",
	title = {Metric Foliations and Curvature},
	publisher = {Birkhäuser Verlag, Basel},
	year = {2009},
}

@Article{R,
	author = {Roitberg, J.},
	title = {On a Construction of Bredon},
	journal = {Proceedings of the American Mathematical Society},
	year = {1972},
	volume = {33},
	pages = {623-626},
}

@PhdThesis{S,
	author = {Sperança, L.},
	title = {Geometria e Topologia de Cobordos},
	school = {Universidade Estadual de Campinas},
	year = {2012},
}

@article{speranca_oddbundles,
  title = {{On Riemannian Foliations over Positively Curved Manifolds}},
  volume = {28},
  ISSN = {1559-002X},
  url = {http://dx.doi.org/10.1007/s12220-017-9901-5},
  DOI = {10.1007/s12220-017-9901-5},
  number = {3},
  journal = {The Journal of Geometric Analysis},
  publisher = {Springer Science and Business Media LLC},
  author = {Speran\c{c}a,  L. D.},
  year = {2017},
  month = jul,
  pages = {2206–2224}
}

@article{wilking2002manifolds,
  title={Manifolds with positive sectional curvature almost everywhere},
  author={Wilking, B.},
  journal={Inventiones mathematicae},
  volume={148},
  number={1},
  pages={117--141},
  year={2002},
  publisher={Springer}
}

\end{document}